\documentclass[10pt]{amsart}
\usepackage{tipa}
\usepackage{amssymb}
\usepackage{stmaryrd}
\usepackage{amsmath,amsfonts,amsthm,mathrsfs,bbm,array,subfigure}

\newtheorem{theorem}{Theorem}[section]

\theoremstyle{definition}

\newtheorem{question}[theorem]{Question}

\theoremstyle{remark}

\numberwithin{equation}{section}

\begin{document}

\title[On the Diophantine equations $z^2=f(x)^2 \pm g(y)^2$]{On the Diophantine equations $z^2=f(x)^2 \pm g(y)^2$ concerning Laurent polynomials}

\author{Yong Zhang}
\address{School of Mathematics and Statistics, Changsha University of Science and Technology,
Changsha 410114, People's Republic of China}
 \email{zhangyongzju$@$163.com}

\thanks{This research was supported by the National Natural Science Foundation of China (Grant No.~11501052).}

\author{Arman Shamsi Zargar}
\address{Young Researchers and Elite Club, Ardabil Branch, Islamic Azad University, Ardabil, Iran}
 \email{shzargar.arman@gmail.com}

%    General info

\date{May 22, 2017}


\begin{abstract}
By the theory of elliptic curves, we study the nontrivial rational
parametric solutions and rational solutions of the Diophantine
equations $z^2=f(x)^2 \pm g(y)^2$ for some simple Laurent
polynomials $f$ and $g$.
\end{abstract}

\maketitle

\section{Introduction}

In 2010, A. Togbe and M. Ulas \cite{Togbe-Ulas} considered the
rational solutions of the Diophantine equations
\begin{equation}\label{Eq11}
z^2=f(x)^2\pm f(y)^2,
\end{equation}
where $f$ being quadratic and cubic polynomials. At the same year,
B. He, A. Togbe and M. Ulas \cite{He-Togbe-Ulas} further
investigated the integer solutions of Eq. (\ref{Eq11}) for some
special polynomials $f$.

In 2016, Y. Zhang \cite{Zhang2016} studied the nontrivial rational
parametric solutions of the Diophantine equations
\[f(x)f(y)=f(z)^n,\]
where $n=1,2,$ concerning the Laurent polynomials $f$. Moreover, we
\cite{Zhang-Zargar} considered the nontrivial rational parametric
solutions and rational solutions of Eq. (\ref{Eq11}) for some simple
Laurent polynomials $f$.

In 2017, Sz. Tengely and M. Ulas \cite{Tengely-Ulas} investigated
the integer solutions of the Diophantine equations
\begin{equation}\label{Eq12}
z^2=f(x)^2+g(y)^2
\end{equation}
and \begin{equation}\label{Eq13} z^2=f(x)^2-g(y)^2
\end{equation}
for different polynomials $f$ and $g$.

In this paper, we continue the study of \cite{Zhang2016} and
\cite{Zhang-Zargar}, and consider the nontrivial rational parametric
solutions and rational solutions of Eqs. (\ref{Eq12}) and
(\ref{Eq13}) for some simple Laurent polynomials $f$ and $g$. The
nontrivial solution $(x,y,z)$ of Eqs. (\ref{Eq12}) and (\ref{Eq13})
respectively means that $f(x)g(y)\neq0$ and $f(x)^2\neq g(y)^2$.

Recall that a Laurent polynomial with coefficients in a field
$\mathbb{F}$ is an expression of the form
\[f(x)=\sum_{k}a_kx^k,a_k\in \mathbb{F},\]
where $x$ is a formal variable, the summation index $k$ is an
integer (not necessarily positive) and only finitely many
coefficients $a_k$ are nonzero. Here we mainly care about the simple
Laurent polynomials
\[f(x)=ax+b+\frac{c}{x},~ax^2+bx+c+\frac{d}{x},~ax+b+\frac{c}{x}+\frac{d}{x^2}\]
with nonzero integers $a,b,c,d$. In the sequel, without losing the
generality we may assume that $a=1$, guaranteed by the shape of Eqs.
(\ref{Eq12}) and (\ref{Eq13}).

By the theory of elliptic curves, we give a positive answer to the
\textbf{Question 3.3} of \cite{Zhang-Zargar}, and prove

\begin{theorem}
For $f=x+a+b/x$ and $g=y^2+cy+d+e/y$ with nonzero integers
$a,b,c,d,e$, Eq. (\ref{Eq12}) or Eq. (\ref{Eq13}) has infinitely
many nontrivial rational parametric solutions.
\end{theorem}

\begin{theorem}
For $f=x+a+b/x$ and $g=y+c+d/y+e/y^2$ with nonzero integers
$a,b,c,d,e$, Eq. (\ref{Eq12}) or Eq. (\ref{Eq13}) has infinitely
many nontrivial rational parametric solutions.
\end{theorem}

\begin{theorem}
For $f=(x+a)(x+b)(x+c)/x$ and $g=(y+d)(y+e)(y+cd/a)/y$ with nonzero
integers $a,b,c,d,e$, Eq. (\ref{Eq12}) or Eq. (\ref{Eq13}) has
infinitely many nontrivial rational solutions.
\end{theorem}

\begin{theorem}
For $f=(x+a)(x+b)(x+c)/x$ and $g=(y+d)(y+e)(y+cd/a)/y^2$ with
nonzero integers $a,b,c,d,e$, Eq. (\ref{Eq12}) or Eq. (\ref{Eq13})
has infinitely many nontrivial rational solutions.
\end{theorem}

\section{The proofs of Theorems}

\begin{proof}[\textbf{Proof of Theorem 1.1.}]
For $f=x+a+b/x$ and $g=y^2+cy+d+e/y$, let \[x=tT,y=T.\] Then Eq.
(\ref{Eq12}) equals
\[\begin{split}
T^2t^2z^2=&T^4t^4+2T^3at^3+[T^6+2cT^5+(c^2+2d)T^4+(2cd+2e)T^3\\
          &+(a^2+2ce+d^2+2b)T^2+2Tde+e^2]t^2+2Tabt+b^2=:g_1(t).\end{split}\]

Let $v=Ttz,$ and consider the curve $\mathcal{C}_1:~v^2=g_1(t)$. In
order to prove Theorem 1.1 we must show that the curve
$\mathcal{C}_1$ has infinitely many $Q(T)$-rational points. The
curve $\mathcal{C}_1$ is a quartic curve with rational point
$P=(0,b)$. By the method of Fermat \cite[p. 639]{Dickson2}, using
the point $P$ we can produce another point $P'=(t_1,v_1)$, which
satisfies the condition $t_1v_1\neq0$. In order to construct a such
point $P'$, we put
\[v= pt^2+ qt +b,\]
where $p,q$ are indeterminate variables. Then
\[v^2-g_1(t)=\sum_{i=1}^4A_it^i,\]
where the quantities $A_i = A_i(p,q)$ are given by
\[\begin{split}
A_1=&2abT-2bq,\\
A_2=&T^6+2cT^5+(c^2+2d)T^4+(2cd+2e)T^3+(a^2+2ce+d^2+2b)T^2\\
&+2deT-2pb+e^2-q^2,\\
A_3=&2aT^3-2pq,\\
A_4=&T^4-p^2.
\end{split}\]
The system of equations $A_4=A_3=0$ in $p, q$ has a solution given
by\[p=-T^2,q=-aT.\]This implies that the equation
\[v^2-g_1(t)=\sum_{i=1}^4A_it^i=0\] has the rational roots $t=0$ and
\[t=-\frac{4abT}{T^6+2cT^5+(c^2+2d)T^4+(2cd+2e)T^3+(2ce+d^2+4b)T^2+2deT+e^2}.\]
Then we have the point $P'=(t_1,v_1)$ with
\[\begin{split}
t_1=&-\frac{4abT}{T^6+2T^5c+(c^2+2d)T^4+(2cd+2e)T^3+(2ce+d^2+4b)T^2+2Tde+e^2},\\
v_1=&\frac{\varphi_1}{\psi_1},
\end{split}\]
where
\[\begin{split}
\varphi_1=&b(T^{12}+4cT^{11}+(6c^2+4d)T^{10}+(4c^3+12cd+4e)T^9+(c^4+12c^2d+4a^2\\
          &+12ce+6d^2+8b)T^8+(4c^3d+8a^2c+12c^2e+12cd^2+16bc+12de)T^7\\
          &+(4a^2c^2+4c^3e+6c^2d^2+8a^2d+8bc^2+24cde+4d^3+16bd+6e^2)T^6\\
          &+(8a^2cd+12c^2de+4cd^3+8a^2e+16bcd+12ce^2+12d^2e+16be)T^5+(8a^2ce\\
          &+4a^2d^2+6c^2e^2+12cd^2e+d^4+16bce+8bd^2+12de^2+16b^2)T^4+4e(2a^2d\\
          &+3cde+d^3+4bd+e^2)T^3+2e^2(2a^2+2ce+3d^2+4b)T^2+4Tde^3+e^4),\\
\psi_1=&(T^6+2cT^5+(c^2+2d)T^4+(2cd+2e)T^3+(2ce+d^2+4b)T^2+2deT+e^2)^{2}.
\end{split}\]

Again, using the method of Fermat, we let \[v= p(t+t_1)^2+ q(t+t_1)
+v_1,\] and get another $P''=(t_2,v_2)$ on the curve
$\mathcal{C}_1$. Repeating the above process any numbers of times
completes the proof of Theorem 1.1 for Eq. (\ref{Eq12}). The same
method can be used to give a proof for Eq. (\ref{Eq13}).

\end{proof}

\begin{proof}[\textbf{Proof of Theorem 1.2.}]
For $f=x+a+b/x$ and $g=y+c+d/y+e/y^2$, let \[x=tT,y=T.\] Then Eq.
(\ref{Eq12}) leads to
\[\begin{split}
T^4t^2z^2=&T^6t^4+2T^5at^3+[T^6+2T^5c+(a^2+c^2+2b+2d)T^4+(2cd+2e)T^3\\
          &+(2ce+d^2)T^2+2Tde+e^2]t^2+2T^3abt+b^2T^2=:g_2(t).\end{split}\]

Let $v=T^2tz,$ and consider the curve $\mathcal{C}_2:~v^2=g_2(t)$.
In order to prove Theorem 1.2 we must show that the curve
$\mathcal{C}_2$ has infinitely many $Q(T)$-rational points. The
curve $\mathcal{C}_2$ is a quartic curve with rational point
$Q=(0,bT)$. By the method of Fermat \cite[p. 639]{Dickson2}, and
using the point $Q$ we can produce another point $Q'=(t_1,v_1)$
satisfying the condition $t_1v_1\neq0$. In order to construct a such
point $Q'$, we take
\[v= pt^2+ qt +bT,\]
where $p,q$ are indeterminate variables. Then
\[v^2-g_2(t)=\sum_{i=1}^4A_it^i,\]
where the quantities $A_i = A_i(p,q)$ are given by
\[\begin{split}
A_1=&2T^3ab-2Tbq,\\
A_2=&T^6+2cT^5+(a^2+c^2+2b+2d)T^4+(2cd+2e)T^3+(2ce+d^2)T^2\\
    &+(2de-2pb)T+e^2-q^2,\\
A_3=&2T^5a-2pq,\\
A_4=&T^6-p^2.
\end{split}\]
The system of equations $A_4=A_3=0$ in $p, q$ has a solution given
by\[p=-T^3,q=-aT^2.\]This implies that the equation
\[v^2-g_2(t)=\sum_{i=1}^4A_it^i=0\] has the rational roots $t=0$ and
\[t=-\frac{4abT}{T^6+2cT^5+(c^2+2d)T^4+(2cd+2e)T^3+(2ce+d^2+4b)T^2+2deT+e^2}.\]
Then we have the point $P'=(t_1,v_1)$ with
\[\begin{split}
t_1=&-\frac{4abT^3}{T^6+2cT^5+(c^2+4b+2d)T^4+(2cd+2e)T^3+(2ce+d^2)T^2+2Tde+e^2},\\
v_1=&\frac{\varphi_1}{\psi_1},
\end{split}\]
where
\[\begin{split}
\varphi_1=&bT(T^{12}+4cT^{11}+(4a^2+6c^2+8b+4d)T^{10}+(8a^2c+4c^3+16bc+12cd\\
&+4e)T^9+(4a^2c^2+c^4+8a^2d+8bc^2+12c^2d+16b^2+16bd+12ce+6d^2)T^8\\
&+(8a^2cd+4c^3d+8a^2e+16bcd+12c^2e+12cd^2+16be+12de)T^7\\
&+(8a^2ce+4a^2d^2+4c^3e+6c^2d^2+16bce+8bd^2+24cde+4d^3+6e^2)T^6\\
&+(8a^2de+12c^2de+4cd^3+16bde+12ce^2+12d^2e)T^5+(4a^2e^2+6c^2e^2\\
&+12cd^2e+d^4+8be^2+12de^2)T^4+4e(3cde+d^3+e^2)T^3\\
&+2e^2(2ce+3d^2)T^2+4Tde^3+e^4),\\
\psi_1=&(T^6+2cT^5+(c^2+4b+2d)T^4+(2cd+2e)T^3+(2ce+d^2)T^2+2Tde+e^2)^{2}.
\end{split}\]

Applying the method of Fermat, we take \[v= p(t+t_1)^2+ q(t+t_1)
+v_1,\] and obtain another point $Q''=(t_2,v_2)$ on the curve
$\mathcal{C}_2$. Repeating the above process any numbers of times
completes the proof of Theorem 1.2 for Eq. (\ref{Eq12}). The same
method can be used to give a proof for Eq. (\ref{Eq13}).
\end{proof}

\begin{proof}[\textbf{Proof of Theorem 1.3.}]
For $f=(x+a)(x+b)(x+c)/x$ and $g=(y+d)(y+e)(y+cd/a)/y$, set
\[x=T,y=\frac{d}{a}T.\] Then Eq. (\ref{Eq12}) reduces to
\[
a^4T^2z^2=(T+a)^2(T+c)^2((a^4+d^4)T^2+(2a^4b+2ad^3e)T+a^2(a^2b^2+d^2e^2)).
\]
To prove the theorem, let us consider the curve
\[\mathcal{C}_3:~(a^4+d^4)T^2+(2a^4b+2ad^3e)T+a^2(a^2b^2+d^2e^2)=S^2.\]
Note that the point $(T,S)=(-b,ade-bd^2)$ lies on $\mathcal{C}_3$,
then it can be parameterized by
\[\begin{split}
T=&\frac{br^2+2d(ae-bd)r-a^4b-2ad^3e+bd^4}{a^4+d^4-r^2},\\
S=&\frac{d(r^2-2d^2r+a^4+d^4)(ae-bd)}{a^4+d^4-r^2},
\end{split}\]
where $r$ is a rational parameter. Hence,
\[\begin{split}
x=&\frac{br^2+2d(ae-bd)r-a^4b-2ad^3e+bd^4}{a^4+d^4-r^2},\\
y=&\frac{d(br^2+2d(ae-bd)r-a^4b-2ad^3e+bd^4)}{a(a^4+d^4-r^2)},\\
z=&\frac{z_1}{z_2},
\end{split}\]
where
\[\begin{split}
z_1=&d(ae-bd)(a^4+d^4-2d^2r+r^2)(a^4b-a^4c+2ad^3e-bd^4-cd^4-2ader+2bd^2r\\
    &+(c-b)r^2)(a^5-a^4b+ad^4-2ad^3e+bd^4+2ader-2bd^2r+(b-a)r^2),\\
z_2=&a^2(a^4+d^4-r^2)^2(a^4b+2ad^3e-bd^4-2ader+2bd^2r-br^2).
\end{split}\]
This completes the proof of Theorem 1.3 for Eq. (\ref{Eq12}). A
similar method for Eq. (\ref{Eq13}).
\end{proof}

\begin{proof}[\textbf{Proof of Theorem 1.4.}]
For $f=(x+a)(x+b)(x+c)/x$ and $g=(y+d)(y+e)(y+cd/a)/y^2$, put
\[x=T,y=\frac{d}{a}T.\] Then Eq. (\ref{Eq12}) equals
\[\begin{split}
T^4a^2z^2=&(T+a)^2(T+c)^2(a^2T^4+2a^2bT^3+(a^2b^2+d^2)T^2+2adeT+a^2e^2).
\end{split}\]
To prove the theorem, consider the curve
\[\mathcal{C}_4:~a^2T^4+2a^2bT^3+(a^2b^2+d^2)T^2+2adeT+a^2e^2=S^2.\]

The curve $\mathcal{C}_4$ is a quartic curve with rational point
$R=(0,ae)$. By the method of Fermat \cite[p. 639]{Dickson2}, and
using the point $R$, we can produce another point $R'=(T_1,S_1)$,
which satisfies the condition $T_1S_1\neq0$. In order to construct a
such point $R'$, we put
\[S= pT^2+ qT +ae,\]
where $p,q$ are indeterminate variables. Then
\[\sum_{i=1}^4A_iT^i=0,\]
where the quantities $A_i = A_i(p,q)$ are given by
\[\begin{split}
A_1=&2ade-2aeq,\\
A_2=&a^2b^2-2aep+d^2-q^2,\\
A_3=&2a^2b-2pq,\\
A_4=&a^2-p^2.
\end{split}\]
The system of equations $A_4=A_3=0$ in $p, q$ has a solution given
by\[p=-a,q=-ab.\]This implies the equation \[\sum_{i=1}^4A_iT^i=0\]
has the rational roots $t=0$ and
\[T=-\frac{2ae(ab+d)}{2a^2e+d^2}.\]
Then we have the point $R'=(T_1,S_1)$, where
\[\begin{split}
T_1=&-\frac{2ae(ab+d)}{2a^2e+d^2},\\
S_1=&\frac{ae(4a^4e^2-4a^3bde+2a^2b^2d^2+2abd^3+d^4)}{(2a^2e+d^2)^2}.
\end{split}\]

Using the method of Fermat, we let
\[S= p(T+T_1)^2+q(T+T_1)+S_1,\]
and find another point $R''=(T_2,S_2)$ on the curve $\mathcal{C}_4$.
Repeat the above process any numbers of times, this completes the
proof of Theorem 1.4 for Eq. (\ref{Eq12}). The same method for
giving a proof for Eq. (\ref{Eq13}).
\end{proof}

\section{Some related questions}

We have studied the rational parametric solutions and rational
solutions of Eqs. (\ref{Eq12}) and (\ref{Eq13}) for
$f,g=x+a+b/x,x^2+ax+b+c/x$ and $x+a+b/x+c/x^2$, but we don't get the
same results for other Laurent polynomials.

\begin{question}
For $f=x^3+ax^2+bx+c+d/x$ or $x^2+ax+b+c/x+d/x^2$, do Eqs.
(\ref{Eq12}) and (\ref{Eq13}) have rational solutions? If they have,
are there infinitely many?
\end{question}

Finding the integer solutions of Eqs. (\ref{Eq12}) and (\ref{Eq13})
is also an interesting question.

\begin{question} Does there exist a Laurent polynomial such that Eqs. (\ref{Eq12}) and
(\ref{Eq13}) have infinitely many integer solutions?
\end{question}

\vskip20pt
\bibliographystyle{amsplain}

\end{document}